\documentclass[11pt]{amsart}
\usepackage{amssymb}

\newtheorem{lemma}{Lemma}[section] 

\newtheorem{theorem}[lemma]{Theorem}

\newtheorem{proposition}[lemma]{Proposition}
\newtheorem{corollary}[lemma]{Corollary}
\newtheorem{definition}[lemma]{Definition}

\begin{document} 
\title{On the topology of semi-algebraic functions on closed semi-algebraic sets}
\author{Nicolas Dutertre}
\address{Universit\'e de Provence, Centre de Math\'ematiques et Informatique,
39 rue Joliot-Curie,
13453 Marseille Cedex 13, France.}
\email{dutertre@cmi.univ-mrs.fr}

\thanks{
 Supported by {\em Agence Nationale de la Recherche} (reference ANR-08-JCJC-0118-01)}

\begin{abstract}
We consider a closed semi-algebraic set $X \subset \mathbb{R}^n$ and a $C^2$ semi-algebraic function $f : \mathbb{R}^n \rightarrow \mathbb{R}$ such that $f_{\vert X}$ has a finite number of critical points. We relate the topology of $X$ to the topology of the sets
$\{ f * \alpha \}$,  where  $* \in \{\le,=,\ge \}$ and $\alpha \in \mathbb{R}$,
and the indices of the critical points of $f_{\vert X}$ and $-f_{\vert X}$. We also relate the topology of $X$ to the topology of the links at infinity of the sets
$\{ f * \alpha \}$ and the indices of these critical points. We give applications when $X=\mathbb{R}^n$ and when $f$ is a generic linear function.
\end{abstract}

\maketitle
\markboth{Nicolas Dutertre}{On the topology of semi-algebraic functions on closed semi-algebraic sets }

\section{Introduction}
Let $f : (\mathbb{R}^n,0) \rightarrow (\mathbb{R},0)$ be an analytic function-germ with an isolated critical point at $0$.
The Khimshiashvili formula [Kh] states that:
$$ \chi \big(f^{-1}(\delta) \cap B_\varepsilon^n \big)=1-\hbox{sign}(-\delta)^n  \hbox{deg}_0 \nabla f,$$
where $0 < \vert \delta \vert \ll \varepsilon \ll 1$, $B_\varepsilon^n $ is a closed ball of radius $\varepsilon$ centered at $0$, 
$\nabla f$ is the gradient of $f$ and deg$_0 \nabla f$ is the topological degree of the mapping $\frac{\nabla f}{\vert \nabla f \vert} :S_\varepsilon^{n-1}  \rightarrow
S^{n-1}$.  As a corollary of this formula, one gets (see [Ar] or [Wa]):
$$\chi \big( \{f \le 0 \} \cap S_\varepsilon^{n-1} \big)=1-\hbox{deg}_0 \nabla f ,$$ 
$$ \chi \big( \{f \ge 0 \}  \cap S_\varepsilon^{n-1} \big)=1+(-1)^{n-1}\hbox{deg}_0 \nabla f,$$ and:
$$\chi \big( \{f = 0 \}  \cap S_\varepsilon^{n-1}\big)=2-2\hbox{deg}_0 \nabla f \hbox{ if } n \hbox{ is even}.$$

In [Se], Sekalski gives a global counterpart of Khimshiashvili's formula for a polynomial mapping $f: \mathbb{R}^2 \rightarrow \mathbb{R}$ with a finite number of critical points. He considers the set 
$\Lambda_f =\{\lambda_1,\cdots,\lambda_k \}$ of critical values of $f$ at infinity, where $\lambda_1 < \ldots < \lambda_k$, and its complement $\mathbb{R} \setminus \Lambda_f = \cup_{i=0}^k ]\lambda_i,\lambda_{i+1}[$ where $\lambda_0=-\infty$ and $\lambda_{k+1}=+\infty$. Denoting by $r_\infty(g)$ the numbers of real branches at infinity of a curve $\{g=0\}$ in $\mathbb{R}^2$, he proves that:
$$\hbox{deg}_\infty \nabla f = 1 +\sum_{i=1}^k r_\infty(f-\lambda_i) -\sum_{i=0}^k r_\infty(f-\lambda_i^+),$$
where for $i=0,\ldots,k$, $\lambda_i^+$ is an element of $]\lambda_i,\lambda_{i+1}[$ and deg$_\infty \nabla f$ is the topological degree of the mapping $\frac{\nabla f}{\Vert \nabla f \Vert} : S_R^{n-1}\rightarrow S^{n-1}$, $R \gg 1$. Here a real branch is homeomorphic to a neighborhood of infinity in $\mathbb{R}$ and hence has two connected components.

Our aim is to generalize Sekalski's formula and to establish other similar results. We consider a closed semi-algebraic set $X \subset \mathbb{R}^n$ equipped with a finite semi-algebraic Whitney stratification $(S_\alpha)_{\alpha \in A}$ and a $C^2$ semi-algebraic function $f : \mathbb{R}^n \rightarrow \mathbb{R}$ such that $f_{\vert X}$ has a finite number of critical points $p_1,\ldots,p_l$.
The index of $f_{\vert X}$ at $p_i$ is defined by:
$$\hbox{ind}(f,X,p_i)=1-\chi \big( X \cap \{f=f(p_i)-\delta \} \cap B_\varepsilon^n(p_i) \big),$$
where $0< \delta \ll \varepsilon \ll 1$. In Section 3, Proposition 3.6 and Corollary 3.7, we give relations between the Euler characteristics of the sets $\{ f * \alpha \}$,  where $* \in \{\le,=,\ge \}$ and $\alpha \in \mathbb{R}$, the indices of the critical points of $f_{\vert X}$ and $-f_{\vert X}$ and four numbers $\lambda_{f,\alpha}$, $\lambda_{-f,-\alpha}$, $\mu_{f,\alpha}$ and $\mu_{-f-,\alpha}$. These numbers are defined in terms of the behavior of $f_{\vert X}$ at infinity (Definition 3.5).  Then we consider the following finite subset of $\mathbb{R}$:
$$\displaylines{
\quad \Lambda^*_f = \Big\{ \alpha \in \mathbb{R} \ \vert \ \beta \mapsto \chi \big( \hbox{Lk}^\infty(X \cap \{ f * \beta \}) \big) \hbox{ is not constant } \hfill \cr
\hfill \hbox{ in a neighborhood of } \alpha \Big\}, \quad \cr
}$$
where $* \in \{ \le, = ,\ge \}$ and Lk$^\infty(-)$ denotes the link at infinity. Writing $\Lambda_f^\le=\{b_1,\ldots,b_r\}$, $\mathbb{R} \setminus \Lambda_f^\le = \cup_{i=0}^r ]b_i,b_{i+1}[$ with $b_0=-\infty$ and $b_{r+1}=+\infty$ and studying the behavior at infinity of the numbers $\lambda_{f,\alpha}$,  $\lambda_{-f,-\alpha}$, $\mu_{f,\alpha}$ and $\mu_{-f-,\alpha}$, we show that  (Theorem 3.16):
$$\displaylines{
\quad \chi(X)=\sum_{i=1}^k \hbox{ind}(f,X,p_i) +\sum_{j=0}^r \chi \big( \hbox{Lk}^\infty(X \cap \{ f \le b_j^+ \}) \big) \hfill \cr
\hfill  - \sum_{j=1}^r \chi \big( \hbox{Lk}^\infty(X \cap \{ f \le b_j \}) \big), \quad \cr
}$$
where for $j \in \{0,\ldots,r\}$, $b_j^+ \in ]b_j,b_{j+1}[$. Similar formulas involving $\Lambda_f^\ge$ and $\Lambda_f^=$ are proved in Theorem 3.17 and Corollary 3.18. 

Next we consider the finite subset $\tilde{B}_f=f\big( \{p_1,\ldots,p_l \} \big) \cup \Lambda_f^\le \cup \Lambda_f^\ge$ of $\mathbb{R}$. We show that if $\alpha \notin \tilde{B}_f$ then the functions $\beta \mapsto \chi \big( X \cap \{f * \beta \} \big)$, $* \in \{\le,=,\ge \}$, are constant in a neighborhood of $\alpha$ (Proposition 3.19). In Theorem 3.20 and Theorem 3.21 we express $\chi(X)$ in terms of the indices of the critical points of $f_{\vert X}$ and $-f_{\vert X}$ and the variations of the Euler characteristics of the sets $\{f * \alpha \}$, $\alpha \in \{\le,=,\ge \}$.   

In Section 4, we apply all these results to the case $X=\mathbb{R}^n$ in order to recover and generalize Sekalski's formula (Theorem 4.4). 

In Section 5, we apply the results of Section 3 to generic linear functions. For $v \in S^{n-1}$, we denote by $v^*$ the function $v^*(x)=\langle v,x \rangle$. We show that for $v$ generic in $S^{n-1}$, the sets $\Lambda_{v^*}^\le$, $\Lambda_{v^*}^=$ and $\Lambda_{v^*}^\ge$ are empty. Hence for such a $v$, the Euler characteristics of the sets $X \cap \{v^* ? \alpha\}$, $? \in \{\le,=,\ge\}$ and $\alpha \in \mathbb{R}$, as well as the Euler characteristics of their links at infinity, can be expressed only in terms of the critical points of $v^*_{\vert X}$ and $-v^*_{\vert X}$ (Proposition 5.4 and Proposition 5.5). We use these results to give a new proof of the Gauss-Bonnet formula for closed semi-algebraic set (Theorem 5.8), that we initially proved in [Dut3, Corollary 5.7] using the technology of the normal cycle [Fu] and a deep theorem of Fu and McCrory [FM, Theorem 3.7]. 

Section 2 of this paper contains three technical lemmas. 

We will use the following notations: for $p \in \mathbb{R}^n$ and $\varepsilon >0$, $B_\varepsilon^n(p)$ is the ball of radius $\varepsilon$ centered at $p$ and $S_\varepsilon^{n-1}(p)$ the sphere of radius $\varepsilon$ centered at $p$. If $p=0$, we simply set $B_\varepsilon^n $ and $S_\varepsilon^{n-1}$ and if $p=0$ and $\varepsilon=1$ we use the standard notations $B^n$ and $S^{n-1}$. If $\mathcal{E}$ is a subset of $\mathbb{R}^n$ then $\mathring{\mathcal{E}}$ denotes its topological interior.

\section{Some lemmas}
Let $X \subset \mathbb{R}^n$ be a closed semi-algebraic set equipped with a semi-algebraic Whitney stratification $(S_\alpha)_{\alpha \in A}$: $X= \sqcup_{\alpha \in A} S_\alpha$. 
Let $g : \mathbb{R}^n \rightarrow \mathbb{R}$ be a $C^2$ semi-algebraic function such that $g^{-1}(0)$ intersects $X$ transversally (in the stratified sense). Then the following partition:
$$X\cap \{ g \le 0 \} = \bigsqcup_{\alpha \in A}  S_\alpha \cap \{g<0\} \sqcup \bigsqcup_{\alpha \in A} S_\alpha \cap \{g=0\},$$
is a Whitney stratification of the closed semi-algebraic set $X \cap \{g \le 0 \}$. 

Let $f : \mathbb{R}^n \rightarrow \mathbb{R}$ be another $C^2$ semi-algebraic function such that $f_{\vert X \cap \{g \le 0 \} }$ admits an isolated critical point $p$ in $X \cap \{g=0\}$ which is not a critical point of $f_{\vert X}$.  If $S$ denotes the stratum of $X$ that contains $p$, this implies that:
$$\nabla (f_{\vert S})(p) =\lambda (p) \nabla (g_{\vert S})(p),$$
with $\lambda(p) \not= 0$. We assume that $f_{\vert \{g=0\}}$ is a submersion in the neighborhood of $p$ if dim $S <n$ and, for simplicity,  that $f(p)=0$. 
%For $\varepsilon >0$, let $B_\varepsilon^n    (p)$ (resp. $S_\varepsilon^{n-1}    (p)$) be the closed ball (resp. the sphere) centered at $p$ with radius $\varepsilon$.

\begin{lemma}
For $0< \delta \ll \varepsilon \ll 1$, we have:
$$\displaylines{
\qquad \chi \big(f^{-1}(-\delta) \cap B_\varepsilon^n    (p) \cap X \cap \{g \le 0 \} \big) =1 \hbox{ if } \lambda(p) >0,  \hfill \cr
\qquad \chi \big(f^{-1}(-\delta) \cap B_\varepsilon^n    (p) \cap X \cap \{g \le 0 \} \big) =  \hfill \cr
 \hfill \chi \big(f^{-1}(-\delta) \cap B_\varepsilon^n    (p) \cap X \cap \{g = 0 \} \big )  \hbox{ if } \lambda(p) <0. \qquad \cr
 }$$
\end{lemma}
{\it Proof.} We assume first that dim $X <n$.  Let $h : \mathbb{R}^n \rightarrow \mathbb{R}$ be a semi-algebraic approximating function for $X$ from outside (see [BK, Definition 6.1]). This implies that $h$ satisfies the following conditions:
\begin{itemize}
\item[(i)] The function $h$ is nonnegative and $h^{-1}\{0\}=X$.
\item[(ii)] The function $h$ is of class $C^3$ on $\mathbb{R}^n \setminus X$.
\item[(iii)] There exists $\delta >0$ such that all $t \in ]0,\delta]$ are regular values of $f$.
\item[(iv)] If $(y_k)_{k \in \mathbb{N}}$ is a sequence of points in $\mathbb{R}^n$ tending to a point $x$ in $X$ such that
$h(y_k) \in ]0,\delta]$ and $\frac{\nabla h}{\Vert \nabla h \Vert}(y_k)$ tends to $v$, then $v$ is normal to $T_x S$, $S$ being the stratum containing $x$.
\end{itemize}
Let us choose $\varepsilon$ sufficiently small so that the ball $B_{\varepsilon'}^n(p)$ intersect $X$, $\{ g \le 0\} $ and $X \cap \{g \le 0 \}$ transversally for $\varepsilon'  \le  \varepsilon$. For $r>0$ sufficiently small, the set $W_{\varepsilon,r} =B_\varepsilon^n    (p) \cap \{g \le 0 \} \cap \{h \le r \}$ is a manifold with corners. To see this, it is enough to prove that $r$ is not a critical value of 
$h_{\vert B_\varepsilon^n    (p) \cap \{ g \le 0 \}}$, which means that $r$ is not a critical value of:
 $$h_1= h_{\vert \mathring{B_\varepsilon^n    (p)} \cap \{ g < 0 \}} \ , \ h_2=h_{\vert \mathring{B_\varepsilon^n    (p)} \cap \{ g =0 \}}\ , \ h_3=h_{\vert S_\varepsilon^{n-1}    (p) \cap \{ g< 0 \}},$$
 and: $$h_4=h_{\vert S_\varepsilon^{n-1}    (p) \cap \{ g = 0 \}}.$$ The fact that $r$ is not a critical value of $h_1$ is trivial by Condition (iii) above. If for $r>0$ small, $r$ is a critical value of $h_2$ then we can find a sequence of points $(q_n)$ in $B_\varepsilon^n    (p) \cap \{g=0\}$ such that $h(q_n) \rightarrow 0$ and $h_{\vert \{g=0\}}$ admits a critical point at $q_n$. Applying Condition (iv) above, we see that there exists a point $q$ in $B_\varepsilon^n    (p) \cap X \cap \{g=0\}$ such that $g^{-1}(0)$ does not intersect $X$ transversally at $q$, which is impossible. Similarly, we can prove that $r$ is not a critical value of $h_3$ and $h_4$.

Let $\delta $ be such that $0 < \delta \ll \varepsilon$ and the fibres $f^{-1}(-\delta)$ and $f^{-1}(\delta)$ intersect $X \cap \{g \le 0 \} \cap B_\varepsilon^n    (p)$ transversally. This is possible since $f$ has an isolated critical point at $p$ on $X \cap \{ g \le 0 \}$. Let us study the critical points of $f_{\vert W_{\varepsilon,r}}$ and $f_{\vert W_{\varepsilon,r} \cap \{g=0\}}$ lying in $f^{-1}([-\delta,\delta])$, when $r$ is small. Always using Condition (iv) above, we can see that they only appear in $\{h=r\} \cap \{g=0\} \cap \mathring{B_\varepsilon^n    (p)}$. Furthermore, with the terminology introduced in [Dut2,\S 2], if $\lambda(p)>0$ then they are all outwards for $f_{\vert W_{\varepsilon,r}}$. If $\lambda (p)<0$ then such a critical point is inwards for $f_{\vert W_{\varepsilon,r}}$ if and only if it is inwards for $f_{\vert W_{\varepsilon,r} \cap \{g=0\}}$. 
Moving $f$ a little, we can assume that these critical points are non-degenerate for $f_{\vert 
\{h=r\} \cap \{g=0\} \cap \mathring{B_\varepsilon^n    (p)}}$. Applying Morse theory for manifolds with corners,  if $\lambda(p)>0$ then we get:
$$\chi \big( f^{-1}([-\delta, \delta]) \cap W_{\varepsilon,r} \big) -\chi \big(f^{-1}(-\delta) \cap W_{\varepsilon,r}\big)=0.$$
If $\lambda(p)<0$ then we get:
$$\chi \big( f^{-1}([-\delta, \delta]) \cap W_{\varepsilon,r} \big) -\chi \big(f^{-1}(-\delta) \cap W_{\varepsilon,r}\big)=$$
$$\chi \big( f^{-1}([-\delta, \delta]) \cap \{g=0\} \cap W_{\varepsilon,r} \big) -\chi \big(f^{-1}(-\delta) \cap \{g=0 \} \cap W_{\varepsilon,r} \big) .$$ 
We conclude remarking that: 
$$\displaylines{
\quad \chi \big(f^{-1}([-\delta,\delta]) \cap W_{\varepsilon,r} \big)= \chi \big(f^{-1}([-\delta,\delta]) \cap X \cap \{g \le 0 \}\cap B_\varepsilon^n    (p) \big)= \hfill \cr
\hfill\chi \big(f^{-1}(0) \cap X \cap \{g \le 0 \}\cap B_\varepsilon^n    (p) \big) =1, \quad \cr
}$$ 
$$\chi \big(f^{-1}(-\delta) \cap W_{\varepsilon,r} \big)= \chi \big(f^{-1}(-\delta) \cap X \cap \{g \le 0 \} \cap B_\varepsilon^n    (p)\big),$$ 
$$\chi \big(f^{-1}([-\delta,\delta]) \cap \{g=0 \} \cap W_{\varepsilon,r} \big)= \chi \big(f^{-1}([-\delta,\delta]) \cap X \cap \{g= 0 \} \cap B_\varepsilon^n    (p)\big)=1,$$ and:
$$\chi \big(f^{-1}(-\delta) \cap \{g=0\} \cap W_{\varepsilon,r} \big)= \chi \big(f^{-1}(-\delta) \cap X \cap \{g =0 \}\cap B_\varepsilon^n    (p) \big),$$ if $r$ is sufficiently small.

If dim $X=n$ then we apply the previous case to the semi-algebraic set $X \times \{ 0 \} \subset \mathbb{R}^{n+1}$ and the functions $F$ and $G$ defined by
$F(x,t)=f(x)+t$ and $G(x,t)=g(x)$, where $(x,t)$ is a coordinate system of $\mathbb{R}^{n+1}=\mathbb{R}^n \times \mathbb{R}$.  

$\hfill \Box$

This lemma was inspired by results on indices of vector fields or 1-forms on stratified sets with boundary (see [KT] or [Sc, Chapter 5]).

Let $M \subset \mathbb{R}^n$ be a $C^2$ semi-algebraic manifold of dimension $k$. Let $f : \mathbb{R}^n \rightarrow \mathbb{R}$ be a $C^2$ semi-algebraic function. Let $\Sigma_f^M$ be the critical set of $f_{\vert M}$. For any $a=(a_1,\ldots,a_n) \in \mathbb{R}^n$, we denote by $\rho_a$ the function $\rho_a(x)=\frac{1}{2} \sum_{i=1}^n (x_i-a_i)^2$ and by $\Gamma^M_{f,a}$ the following semi-algebraic set:
$$\Gamma^M_{f,a}=\left\{ x \in M \ \vert \ \hbox{rank} \big[ \nabla(f_{\vert M})(x), \nabla({\rho_a}_{\vert M})(x) \big] < 2 \right\}.$$
\begin{lemma}
For almost all $a \in \mathbb{R}^n$, $\Gamma^M_{f,a} \setminus \Sigma_f^M$ is a smooth semi-algebraic curve (or empty).
\end{lemma}
{\it Proof.} Let $Z$ be the semi-algebraic set of $\mathbb{R}^n \times \mathbb{R}^n$ defined as follows:
$$\displaylines{
\qquad Z =\Big\{ (x,a) \in \mathbb{R}^n \times \mathbb{R}^n \ \vert \  x \in M\setminus \Sigma_f^M \hbox{ and } \hfill \cr
\hfill  \hbox{ rank}\big[ \nabla(f_{\vert M})(x), \nabla({\rho_a}_{\vert M})(x) \big] < 2 \Big\}. \qquad \cr
}$$
Let $(x,a)$ be a point in $Z$. We can suppose that around $x$, $M$ is defined by the vanishing of $l=n-k$ semi-algebraic functions $g_1,\ldots,g_l$ of class $C^2$. Hence in a neighborhood of $(x,a)$, $Z$ is defined by the vanishing of $g_1,\ldots,g_l$ and the minors:
$$\frac{\partial(g_1,\ldots,g_l,f,\rho_a)}{\partial(x_{i_1},\ldots,x_{i_{l+2}})}.$$
Furthermore since $x $ belongs to $M \setminus \Sigma_f^M$, we can assume that:
$$\frac{\partial(g_1,\ldots,g_l,f)}{\partial(x_1,\ldots,x_l,x_{l+1})}(x) \not= 0.$$
Therefore $Z$ is locally defined by $g_1=\ldots=g_l=0$ and:
$$\frac{\partial(g_1,\ldots,g_l,f,\rho_a)}{\partial(x_1,\ldots,x_{l+1},x_{l+2})}= \dots= \frac{\partial(g_1,\ldots,g_l,f,\rho_a)}{\partial(x_1,\ldots,x_{l+1},x_{n})}=0,$$
(see [Dut1,\S 5] for a proof of this fact). Since the gradient vectors of these functions are linearly independent, we see that $Z$ is a smooth semi-algebraic manifold of dimension $2n-(l+n-(l+2)+1)=n+1$. Now let us consider the projection $\pi_2 : Z \rightarrow \mathbb{R}^n$, $(x,a) \mapsto a$. Bertini-Sard's theorem (see [BCR, Th\'eor\`eme 9.5.2]) implies that the set $D_{\pi_2}$ of critical values of $\pi_2$ is a semi-algebraic set of dimension strictly less than $n$. Hence, for all $a \notin D_{\pi_2}$, $\pi_2^{-1}(a)$ is a smooth semi-algebraic curve (maybe empty). But this set is exactly $\Gamma^M_{f,a} \setminus \Sigma_f^M$.   $\hfill \Box$ 

Now consider a semi-algebraic set $Y \subset M$ of dimension strictly less than $k$. We will need the following lemma.
\begin{lemma}
For almost all $a \in \mathbb{R}^n$, $(\Gamma^M_{f,a} \setminus \Sigma^M_f) \cap Y$ is a semi-algebraic set of dimension at most $0$.
\end{lemma}
{\it Proof.} Since $Y$ admits a finite Whitney semi-algebraic stratification, we can assume that $Y$ is smooth of dimension $d<k$. Let $W$ be the semi-algebraic set of $\mathbb{R}^n \times
\mathbb{R}^n$ defined by:
$$\displaylines{
\qquad W = \Big\{ (x,a) \in \mathbb{R}^n \times \mathbb{R}^n \ \vert \ x \in Y \setminus \Sigma^M_f \hbox{ and } \hfill \cr
\hfill \hbox{ rank}\big[ \nabla (f_{\vert M})(x), \nabla ({\rho_a}_{\vert M})(x) \big] <2 \Big\}. \qquad \cr
}$$
Using the same method as in the previous lemma, we can prove that $W$ is a smooth semi-algebraic manifold of dimension $n+1+d-k$.
%$$2n-(n-d+n-(n-k+2)+1)=2n-(n+k-d-1)=n+1+d-k.$$
We can conclude as in the previous lemma, remarking that $d-k \le -1$. $\hfill \Box$

\section{Topology of semi-algebraic functions}

For any closed semi-algebraic set equipped with a Whitney stratification $X = \sqcup_{\alpha \in A} S_\alpha$, we denote by Lk$^\infty(X)$ the link at infinity of $X$. It is defined as follows. Let $\rho : \mathbb{R}^n \rightarrow \mathbb{R}$ be a $C^2$ proper semi-algebraic positive function. Since $\rho_{\vert X}$ is proper, the set of critical points of $\rho_{\vert X}$ (in the stratified sense) is compact. Hence for $R$ sufficiently big, the map $\rho : X \cap \rho^{-1}([R,+\infty[) \rightarrow \mathbb{R}$ is a stratified submersion. 
The link at infinity of $X$ is the fibre of this submersion. The topological type of Lk$^\infty(X)$ does not depend on the choice of the function $\rho$. Indeed, if $\rho_0$ and $\rho_1 : \mathbb{R}^n \rightarrow \mathbb{R}$ are two $C^2$ proper semi-algebraic functions then $X \cap \rho_0^{-1}(R_0)$ and $X \cap \rho_1^{-1}(R_1)$ are homeomorphic for $R_0$ and $R_1$ big enough. To see this, we can apply the procedure described by Durfee in [Dur]. First we remark that, applying the Curve Selection Lemma at infinity [NZ, Lemma 2], for each stratum $S_\alpha$ of $X$, the gradient vector fields $\nabla({\rho_0}_{\vert S_\alpha})$ and $\nabla({\rho_1}_{\vert S_\alpha})$ do not point in opposite direction in a neighborhood of infinity. Next we choose $R_0$ and $R_1$ sufficiently big so that $\rho_0^{-1}
([R_0,+\infty[) \subset \rho_1^{-1} ([R_1,+\infty[) $ and all the gradient vector fields $\nabla({\rho_0}_{\vert S_\alpha})$ and $\nabla({\rho_1}_{\vert S_\alpha})$ do not point in opposite direction in 
$\rho_0^{-1} ([R_0,+\infty[)$. Then the function $\rho : \big( \rho_0^{-1}([R_0,+\infty[) \setminus \rho_1^{-1} ([R_1,+\infty[) \big) \cap X  \rightarrow [0,1]$ defined by:
$$\rho(x)=\frac{R_1-\rho_1(x)}{R_1-\rho_1(x) + \rho_0(x) -R_0},$$
is a proper stratified submersion such that $\rho^{-1}(0)= X \cap \rho_1^{-1}(R_1)$ and $\rho^{-1}(1)=X \cap \rho_0^{-1}(R_0)$. 

Let $f : \mathbb{R}^n \rightarrow \mathbb{R}$ be a $C^2$ semi-algebraic function such that $f_{\vert X} : X \rightarrow \mathbb{R}$ has a finite number of critical points (in the stratified sense) $p_1,\ldots,p_l$. For each $p_i$, we define the index of $f_{\vert X}$ at $p_i$ as follows:
$$\hbox{ind}(f,X,p_i)=1-\chi \big( X \cap \{f=f(p_i)-\delta \} \cap B_\varepsilon^n    (p_i) \big),$$
where $0< \delta \ll \varepsilon \ll 1$. Since we are in the semi-algebraic setting, this index is well-defined thanks to Hardt's  theorem [Ha]. The following theorem is well-known.
\begin{theorem}
If $f_{\vert X}$ is proper then for any $\alpha \in \mathbb{R}^n$, we have:
$$\chi\big(X \cap \{ f \ge \alpha\} \big)-\chi \big(X \cap \{ f = \alpha \} \big) = \sum_{i : f(p_i) > \alpha} \hbox{\em ind}(f,X,p_i),$$ and:
$$\chi \big(X \cap \{ f \le \alpha\} \big)- \chi \big( \hbox{\em Lk}^\infty(X \cap \{ f \le \alpha \} ) \big)= \sum_{i  : f(p_i) \le  \alpha} \hbox{\em ind}(f,X,p_i).$$
\end{theorem}
{\it Proof.} We use Viro's method of integration with respect to the Euler characteristic with compact support, denoted by $\chi_c$. 

For all $x \in X$, let $\varphi(x)=\chi_c \big(X \cap f^{-1}(x^-) \cap B_\varepsilon^n      (x) \big)$ where $x^-$ is a regular value of $f$ close to $f(x)$ with $x^- \le f(x)$. Applying Fubini's theorem  [Vi, Theorem 3.A] to the restriction of $f$ to $X \cap \{ f > \alpha \}$, we get:
$$\int_{X \cap \{f > \alpha \}} \varphi(x) d \chi_c(x) = \int_{]\alpha,+\infty[ } \left( \int_{f^{-1}(y)} \varphi(x) d\chi_c(x) \right) d\chi_c(y).$$
For any $y \in \mathbb{R}$, let $y^-$ be a regular value of $f_{\vert X}$ close to $y$ with $y^- \le y$. Let us denote by $z_1,\ldots,z_s$ the critical points of $f_{\vert X}$ lying in $f^{-1}(y)$. We have:
$$\displaylines{
\qquad \chi_c \big(X \cap f^{-1}(y^-) \big)=\chi_c \big(X \cap f^{-1}(y^-) \setminus \cup_{i=1}^s B_\varepsilon^n    (z_i) \big)+ \hfill \cr
\hfill  \sum_{i=1}^s \chi_c \big(X \cap f^{-1}(y^-) \cap B_\varepsilon^n    (z_i) \big)= \qquad \cr
\qquad \chi_c \big(X \cap f^{-1}(y) \setminus \cup_{i=1}^s B_\varepsilon^n    (z_i) \big)+\sum_{i=1}^s  \varphi(z_i) = \hfill \cr
\hfill \chi_c \big(X \cap f^{-1}(y) \setminus \{z_1,\ldots,z_s \} \big)+\sum_{i=1}^s \varphi (z_i)= \qquad \cr$$
\qquad  \int_{X \cap f^{-1}(y) \setminus \{z_1,\ldots,z_s \}} \varphi(x) d\chi_c(x) + \sum_{i=1}^s  \varphi (z_i) = \int_{f^{-1}(y)} \varphi(x) d\chi_c(x). \hfill \cr
}$$
Let us write:
$$ ]\alpha, +\infty[=]\alpha,\alpha_1 ] \ \cup \  ]\alpha_1,\alpha_2 ]  \  \cup \ldots \cup \  ]\alpha_{j-1},\alpha_j] \   \cup \ ]\alpha_j,+\infty[,$$
 where $\alpha_1,\ldots,\alpha_j$ are the critical values of $f_{\vert X}$ strictly greater than $\alpha$. Since $\chi_c(]\alpha_k, \alpha_{k+1}])=0$ and $f_{\vert X \cap ]\alpha_k, \alpha_{k+1} [}$ is locally trivial, we obtain that:
$$ \int_{X \cap \{f > \alpha \}} \varphi(x) d \chi_c(x) = -\chi_c (X \cap f^{-1}(\beta)),$$
where $\beta$ is a regular value of $f$ strictly greater than $\alpha_j$ and therefore:
$$\chi_c \big(X \cap \{ f > \alpha \} \big)+\chi_c \big(X \cap f^{-1}(\beta) \big)= \sum_{i : f(p_i) > \alpha} \hbox{ind}(f,X,p_i).$$
Applying this equality to $\alpha=\beta$ and using the local triviality of $f_{\vert X}$ over $[\beta,+\infty[$, we get:
$$\chi_c \big(X \cap \{ f > \beta \} \big) + \chi_c \big( X \cap f^{-1}(\beta) \big)=0.$$
Therefore:
$$\chi_c \big(X \cap \{ f > \alpha \} \big)+\chi_c \big( X\cap f^{-1}(\beta) \big)=$$
$$\chi_c \big(X \cap \{\alpha \le f \le \beta \} \big) -\chi_c \big(X \cap \{ f =\alpha \} \big)- \chi_c \big(X \cap \{ f > \beta \} \big) +\chi_c \big(X \cap f^{-1}(\beta) \big)=$$
$$\chi \big(X \cap \{\alpha \le f \le \beta \} \big)-\chi \big(X \cap \{ f =\alpha \} \big).$$
To conclude, we remark that, since $f_{\vert X \cap \{ f \ge \alpha \}}$ is proper and locally trivial over $[\beta, +\infty[$, $X \cap \{\alpha \le f \le \beta \}$ is a deformation retract of 
$X \cap \{ \alpha \le f \}$.   
The second equality is proved  with the same method and the fact that $\chi_c(Y)= \chi(Y) -\chi \big( \hbox{Lk}^\infty(Y) \big)$ for any closed semi-algebraic set $Y \subset \mathbb{R}^n$. $\hfill \Box$

The following corollaries are straightforward consequences of the previous theorem.

\begin{corollary}
If $f_{\vert X}$ is proper then for any $\alpha \in \mathbb{R}$, we have:
$$\displaylines{
\quad \chi \big(X \cap \{f= \alpha \} \big)=\chi(X)-  
\sum_{i : f(p_i) > \alpha} \hbox{\em ind}(f,X,p_i) - \sum_{i : f(p_i) < \alpha} \hbox{\em ind}(-f,X,p_i) , \qquad \cr
}$$
and:
$$\displaylines{
\quad \chi \big(X \cap \{ f \ge \alpha\} \big) -\chi \big(X \cap \{ f \le \alpha\} \big)=\sum_{i : f(p_i) > \alpha} \hbox{\em ind}(f,X,p_i) - \hfill \cr
\hfill \sum_{i: f(p_i) < \alpha} \hbox{\em ind}(f,X,p_i). \qquad \cr
}$$
\end{corollary}
$\hfill \Box$

\begin{corollary}
If $f_{\vert X}$ is proper then for any $\alpha \in \mathbb{R}^n$, we have :
$$\chi \big( \hbox{\em Lk}^{\infty} (X \cap \{f \le \alpha \} ) \big) =\chi(X)-\sum_{i=1}^l  \hbox{\em ind}(f,X,p_i).$$
\end{corollary}
$\hfill \Box$

\begin{corollary} 
If $f_{\vert X}$ is proper then we have:
$$2 \chi(X)- \chi \big(\hbox{\em Lk}^\infty(X) \big)= \sum_{i=1}^l  \hbox{\em ind}(f,X,p_i)+ \sum_{i=1}^l  \hbox{\em ind}(-f,X,p_i).$$
\end{corollary}
$\hfill \Box$

Now we want to investigate the case when $f_{\vert X}$ is not proper. Keeping the notations of the previous section, for $a \in \mathbb{R}^n$, we define $\Gamma^X_{f,a}$ and $\Gamma_{f,a}$ by:
$$\displaylines{
\qquad \Gamma^X_{f,a}=\Big\{ x \in X \ \vert \ \hbox{rank} \big[ \nabla(f_{\vert S})(x), \nabla({\rho_a}_{\vert S})(x) \big] < 2 \hfill \cr
\hfill  \hbox{ where } S \hbox{ is the stratum that contains } x \Big\}, \qquad \cr
}$$
$$\Gamma_{f,a}=\left\{ x \in \mathbb{R}^n \ \vert \ \hbox{rank} \left[ \nabla f(x), \nabla \rho_a(x) \right] < 2 \right\}.$$
By Lemma 2.2, we can choose $a$ such that $\Gamma^X_{f,a}$ is a smooth semi-algebraic curve outside a compact set of $X$. Applying Lemma 2.3 to $M=\mathbb{R}^n$ and $Y$ the closed semi-algebraic set defined as the union of the strata of $X$ of dimension strictly less than $n$, we can choose $a$ such that $\Gamma^X_{f,a}$ and $\Gamma_{f,a}$ do not intersect outside a compact set of $Y$. 
Let us fix $\alpha \in \mathbb{R}$ and $R \gg 1$ such that:
\begin{enumerate}
\item $X \cap B_R^n(a)$ is a deformation retract of $X$,
\item $X \cap \{ f * \alpha \} \cap B_R^n(a)    $ is a deformation retract of $X \cap \{ f * \alpha \} $ where $* \in \{\le, = , \ge \}$,
\item $S_R^{n-1}    (a)    $ intersects $X$ and $X \cap \{ f * \alpha \}$ transversally,
\item $\Gamma^X_{f,a} \cap S_R^{n-1}    (a)    $ is a finite set of points $q_1^R,\ldots,q_m^R$,
\item $p_1,\ldots,p_l \in \mathring{B_R^n(a)    }.$
\end{enumerate}
For each $j \in \{1,\ldots,m\}$, $q_j^R$ is a critical point of $f_{\vert X \cap S_R^{n-1}    (a)    }$ but not a critical point of $f_{\vert X}$. Hence there exists $\mu_j^R \not= 0$ such that:
$$\nabla (f_{\vert S})(q_j^R)=\mu_j^R \nabla({\rho_a}_{\vert S})(q_j^R),$$
where $S$ is the stratum that contains $q_j^R$.

\begin{definition}
We set:
$$\lambda_{f,\alpha} = \sum_{j : f(q_j^R) > \alpha \atop \mu_j^R <0} \hbox{\em ind} (f,X\cap S_R^{n-1}    (a) ,q_j^R),$$
$$\mu_{f,\alpha} = \sum_{j : f(q_j^R) < \alpha \atop \mu_j^R >0} \hbox{\em ind} (f,X\cap S_R^{n-1}    (a) ,q_j^R),$$
$$\nu_{f,\alpha} = \sum_{j : f(q_j^R) < \alpha \atop \mu_j^R <0} \hbox{\em ind} (f,X\cap S_R^{n-1}    (a)  ,q_j^R).$$
\end{definition}

The fact that $\lambda_{f,\alpha}$, $\mu_{f,\alpha}$ and $\nu_{f,\alpha}$ do not depend on $R$ will appear in the next propositions.
Let us remark that if $\mu_j^R<0$ then $f(q_j^R)$ decreases to $-\infty$ or to a finite value as $R$ tends to $+\infty$, and that if $\mu_j^R>0$ then $f(q_j^R)$ increases to $+\infty$ or to a finite value as $R$ tends to $+\infty$. This implies that when $f_{\vert X}$ is proper, the numbers $\lambda_{f,\alpha}$ and $\mu_{f,\alpha}$ vanish.
\begin{proposition}
For any $\alpha \in \mathbb{R}$, we have:
$$\chi \big(X \cap \{f \ge \alpha \} \big) - \chi \big(X \cap \{f = \alpha \} \big)= \sum_{i : f(p_i) > \alpha} \hbox{\em ind}(f,X,p_i) + \lambda_{f,\alpha},$$
and:
$$\chi \big(X \cap \{f \le \alpha \} \big) - \chi \big(X \cap \{f = \alpha \} \big)= \sum_{i :f(p_i) < \alpha} \hbox{\em ind}(-f,X,p_i) + \lambda_{-f,-\alpha}.$$
\end{proposition}
{\it Proof.} We apply Theorem 3.1 to $f_{\vert X \cap B_R^n(a)}$ and we get:
$$\displaylines{
\qquad \chi \big(X \cap B_R^n(a)     \cap \{ f \ge \alpha \} \big)-\chi \big(X \cap B_R^n(a) \cap \{f=\alpha \} \big)= \hfill \cr
\hfill \sum_{i:f(p_i) > \alpha} \hbox{ind}(f,X,p_i) + \sum_{j : f(q_j^R) > \alpha} \hbox{ind}(f,X\cap B_R^n(a) ,q_j^R), \qquad
}$$
and:
$$\displaylines{
\qquad \chi \big(X \cap B_R^n(a)   \cap \{ f \le \alpha \} \big)-\chi \big(X \cap B_R^n(a)  \cap \{f=\alpha \} \big)= \hfill \cr
\hfill \sum_{i :f(p_i) <\alpha} \hbox{ind}(-f,X,p_i) + \sum_{j : f(q_j^R) < \alpha} \hbox{ind}(-f,X\cap B_R^n(a) ,q_j^R).
}$$
Since $\Gamma^X_{f,a}$ and $\Gamma_{f,a}$ do not intersect outside a compact set of $Y$, we can use Lemma 2.1 to evaluate ind($f,X\cap B_R^n(a) ,q_j^R)$ and ind$(-f,X \cap B_R^n(a) ,q_j^R)$. Namely, if $\mu_j^R >0$ then we have: 
$$\hbox{ind}(f,X \cap B_R^n(a),q_j^R)=0,$$
and: $$ \hbox{ind}(-f,X \cap B_R^n(a),q_j^R)=\hbox{ind}(-f,X \cap S_R^{n-1}    (a),q_j^R).$$
If $\mu_j^R <0$ then we have: 
$$\hbox{ind}(f,X \cap B_R^n(a),q_j^R)= \hbox{ind}(f,X \cap S_R^{n-1}    (a),q_j^R),$$
and:
$$ \hbox{ind}(-f,X \cap B_R^n(a)  ,q_j^R)=0.$$
Moreover, by our choice on $R$, $\chi \big(X \cap B_R^n(a) \cap \{f * \alpha \} \big)=\chi \big(X \cap \{f * \alpha \} \big)$ for $* \in \{\le,=,\ge\}$. $\hfill \Box$

\begin{corollary}
For any $\alpha \in \mathbb{R}$, we have:
$$\displaylines{
\quad \chi \big(X \cap \{f=\alpha \} \big)= \chi(X) - \sum_{i:  f(p_i) > \alpha} \hbox{\em ind}(f,X,p_i) - \sum_{i :f(p_i) < \alpha} \hbox{\em ind}(-f,X,p_i)  \hfill \cr
\hfill -\lambda_{f,\alpha}-\lambda_{-f,-\alpha}, \qquad \cr
}$$
and:
$$\displaylines{
\quad \chi \big( X \cap \{f \ge \alpha \} \big)-\chi \big( X \cap \{f \le \alpha \} \big ) =\sum_{i :f(p_i) > \alpha} \hbox{\em ind}(f,X,p_i) +\lambda_{f,\alpha} \hfill \cr
\hfill -\sum_{i : f(p_i) > \alpha} \hbox{\em ind}(f,X,p_i) -\lambda_{-f,-\alpha}. \qquad \cr
}$$
\end{corollary}
$\hfill \Box$

It is also possible to write indices formulas for $\chi \big(\hbox{Lk}^\infty(X \cap \{ f * \alpha \} \big)$, $* \in \{\le,=,\ge\}$.
\begin{proposition}
For any $\alpha \in \mathbb{R}$, we have:
$$\chi\big(\hbox{\em Lk}^{\infty} (X \cap \{f \le \alpha \}\big)=\chi(X) - \sum_{i=1}^l \hbox{\em ind}(f,X,p_i) -\lambda_{f,\alpha}+\mu_{f,\alpha},$$
$$\chi \big(\hbox{\em Lk}^{\infty} (X \cap \{f \ge \alpha \} \big)=\chi(X) - \sum_{i=1}^l \hbox{\em ind}(-f,X,p_i) -\lambda_{-f,-\alpha}+\mu_{-f,-\alpha}.$$
\end{proposition}
{\it Proof.} By Theorem 3.1 applied to $f_{\vert X \cap S_R^{n-1}(a)}$, we have:
$$\chi \big( \hbox{Lk}^{\infty}(X \cap \{f \le \alpha \} \big)=\mu_{f,\alpha}+ \nu_{f,\alpha},$$ and, by Corollary 3.3  applied to  $f_{\vert X \cap B_R^n(a)}$ and by Lemma 2.1:
$$\chi(X)=\sum_{i=1}^l \hbox{ind}(f,X,p_i) +\lambda_{f,\alpha}+\nu_{f,\alpha},$$
because $f^{-1}(\alpha)$ intersect $X \cap S_R^{n-1} (a)    $ transversally.
Similarly, we can write:
$$\chi\big( \hbox{Lk}^{\infty}(X \cap \{f \ge \alpha \} \big)=\mu_{-f,-\alpha}+ \nu_{-f,-\alpha},$$ and:
$$\chi(X)=\sum_{i=1}^l \hbox{ind}(-f,X,p_i) +\lambda_{-f,-\alpha}+\nu_{-f,-\alpha}.$$
$\hfill \Box$

\begin{corollary}
For any $\alpha \in \mathbb{R}$, we have:
$$\displaylines{
\quad \chi \big(\hbox{\em Lk}^{\infty}(X \cap \{f = \alpha \} \big)=2\chi(X) - \chi \big(\hbox{\em Lk}^{\infty}(X) \big) -\sum_{i=1}^l \hbox{\em ind}(f,X,p_i) \hfill \cr
\hfill -\sum_{i=1}^l \hbox{\em ind}(-f,X,p_i) -\lambda_{f,\alpha}+\mu_{f,\alpha} -\lambda_{-f,-\alpha}+\mu_{-f,-\alpha}. \quad \cr
}$$
\end{corollary} 
$\hfill \Box$

In the sequel, we will use these results to establish relations between $\chi(X)$, the indices of the critical points of $f_{\vert X}$ and $-f_{\vert X}$ and the variations of the Euler characteristics $\chi\big( \hbox{Lk}^\infty (\{f * \alpha \}) \big)$, where $* \in \{\le,=,\ge \}$ and $\alpha \in \mathbb{R}$. 
We start with definitions.
\begin{definition}
Let $\Lambda_f$ be the following set:
$$\Lambda_f =\big\{ \alpha \in \mathbb{R} \ \vert \ \exists (x_n)_{n \in \mathbb{N}} \hbox{ in } \Gamma^X_{f,a} \hbox{ such that } \Vert x_n \Vert \rightarrow + \infty \hbox{ and }  f(x_n) \rightarrow \alpha \big\}.$$
\end{definition}
Since $\Gamma^X_{f,a}$ is a curve, $\Lambda_f$ is clearly a finite set. The set $\Lambda_f$ was introduced and studied by Tibar [Ti1] when $X=\mathbb{R}^n$ and $f:\mathbb{R}^n \rightarrow \mathbb{R}$ is a polynomial. Following his terminology, $\Lambda_f$ is the set of points $\alpha$ such that the fibre $f^{-1}(\alpha)$ is not $\rho_a$-regular.

\begin{definition}
Let $* \in \{ \le, = ,\ge \}$. We define $\Lambda^*_f$ by:
$$\displaylines{
\quad \Lambda^*_f = \Big\{ \alpha \in \mathbb{R} \ \vert \ \beta \mapsto \chi \big( \hbox{\em Lk}^\infty(X \cap \{ f * \beta \}) \big) \hbox{ is not constant } \hfill \cr
\hfill \hbox{ in a neighborhood of } \alpha \Big\}. \quad \cr
}$$
\end{definition}

\begin{lemma} 
The sets $\Lambda^\le_f$, $\Lambda^=_f$ and $\Lambda^\ge_f$ are included in $\Lambda_f$.
\end{lemma}
{\it Proof.} If $\alpha$ does not belong to $\Lambda_f$ then we can find a small interval $]-\delta + \alpha, \alpha + \delta [$ such that $\Gamma^X_{f,a}$ and $f^{-1}(]-\delta + \alpha, \alpha + \delta [) \cap X$ do not intersect outside a compact set of $X$. Then we can choose $R \gg 1$ such that for all $\beta \in ]-\delta + \alpha, \alpha + \delta [$, $\hbox{Lk}^\infty(X \cap \{f * \beta \})= X \cap \{f * \beta \} \cap S_R^{n-1}    (a)    $. But $f$ has no critical point in $X \cap \{ -\delta+\alpha < f < \alpha+\delta \} \cap S_R^{n-1}    (a)$, so the Euler characteristics $\chi(X \cap \{ f * \beta \})$ are constant in $]-\delta + \alpha, \alpha+\delta [$. $\hfill \Box$

\begin{corollary}
The sets $\Lambda_f^\le$, $\Lambda_f^=$ and $\Lambda_f^\ge$ are finite.
\end{corollary}
 $\hfill \Box$
\begin{lemma} We have: 
$\Lambda_f^= \subset \Lambda_f^\le \cup \Lambda_f^\ge$, $\Lambda_f^\le \subset \Lambda_f^= \cup \Lambda_f^\ge$,
$\Lambda_f^\ge \subset \Lambda_f^\le \cup \Lambda_f^=$.
\end{lemma}
{\it Proof.} If $\alpha \notin \Lambda_f^\le \cup \Lambda_f^\ge$, then $\beta \mapsto \chi \big(\hbox{Lk}^\infty (X \cap \{f \le \beta \})\big)$ and $\beta \mapsto \chi \big(\hbox{Lk}^\infty (X \cap \{f \ge \beta \})\big)$ are constant in an interval $]-\delta+\alpha,\alpha + \delta [$. By the Mayer-Vietoris sequence, $\beta \mapsto \chi\big( \hbox{Lk}^\infty(X \cap \{f=\beta \}) \big)$ is also constant in $]-\delta+\alpha,\alpha+\delta [$. $\hfill \Box$

\begin{corollary}
We have:  $\Lambda_f^\le \cup \Lambda_f^\ge = \Lambda_f^\le \cup \Lambda_f^= = \Lambda_f^= \cup \Lambda_f^\ge$.
\end{corollary}
 $\hfill \Box$

Since $\Lambda_f^\le$ is finite, we can write $\Lambda_f^\le = \{b_1,\ldots,b_r\}$ where $b_1 < b_2 < \ldots < b_r$ and:
$$\mathbb{R} \setminus \Lambda_f^\le =]-\infty, b_1[ \  \cup \ ]b_1,b_2[ \  \cup \cdots \cup \  ]b_{r-1},b_r [ \ \cup \  ]b_r, +\infty [.$$
On each connected component of $\mathbb{R} \setminus \Lambda_f^\le$, the function $\beta \mapsto \chi \big(\hbox{Lk}^\infty(X \cap \{f \le \beta \} \big)$ is constant. For each $j \in \{0,\ldots,r\}$, let $b_j^+$ be an element of $]b_j,b_{j+1}[$ where $b_0=-\infty$ and $b_{r+1}= + \infty$. 
\begin{theorem}
We have:
$$\displaylines{
\quad \chi(X)=\sum_{i=1}^k \hbox{\em ind}(f,X,p_i) +\sum_{j=0}^r \chi \big( \hbox{\em Lk}^\infty(X \cap \{ f \le b_j^+ \}) \big) \hfill \cr
\hfill  - \sum_{j=1}^r \chi \big( \hbox{\em Lk}^\infty(X \cap \{ f \le b_j \}) \big). \quad \cr
}$$
\end{theorem} 
{\it Proof.} Assume first that $\Lambda_f=\Lambda_f^\le$. Let us choose $R \gg 1$ such that $X \cap B_R^n(a) $ is a deformation retract of $X$, 
$\{p_1,\ldots,p_l\} \subset \mathring{B_R^n(a)}$ and:
$$\Gamma_{f,a}^X \cap \big(\mathbb{R}^n \setminus \mathring{B_R^n(a)} \big)=\sqcup_{i=j}^m \mathcal{B}_j.$$
We have $\Gamma^X_{f,a} \cap S_R^{n-1}    (a)     = \{q_1^R,\ldots,q_m^R\}$. Let us recall that:  $$\nabla( f_{\vert S}) (q_j^R)=\mu_j^R 
\nabla( {\rho_a}_{\vert S})(q_j^R),$$ where $S$ is the stratum that contains $q_j$. By Corollary 3.3 and Lemma 2.1, we can write:
$$\chi(X)=\sum_{i=1}^l \hbox{ind}(f,X,p_i) + \sum_{j : \mu_j^R <0} \hbox{ind}(f,X \cap S_R^{n-1}    (a),q_j^R).$$
We can decompose the second sum in the right hand side of this equality into: 
$$\sum_{j:\mu_j^R <0\atop f(q_j^R) \rightarrow -\infty} \hbox{ind}(f,X \cap S_R^{n-1}    (a),q_j^R),$$
and:
$$\sum_{i=1}^r \sum_{j : \mu_j^R<0 \atop f(q_j^R) \rightarrow b_i} \hbox{ind}(f,X \cap S_R^{n-1}    (a),q_j^R).$$
Let us fix $i $ in $\{1,\ldots,r \}$ and evaluate $\sum_{j : \mu_j^R<0 \atop f(q_j^R) \rightarrow b_i} \hbox{ind}(f,X \cap S_R^{n-1}    (a),q_j^R)$. Since $\mu_j^R<0$, the points $q_j^R$ lie in $\{ f > b_i \}$. Let us choose $R \gg 1$ and $b_i^+$ close to $b_i$ in $]b_i,b_{i+1}[$ such that:
$$\left[  \cup_{j : f \rightarrow b_i \atop along \ \mathcal{B}_j}  \mathcal{B}_j \right] \cap \left\{ \Vert x-a \Vert \ge R \right\} \subset f^{-1}(]b_i,b_i^+[) \cap \left\{ \Vert x-a \Vert \ge R \right\} ,$$
and $X \cap \{ f \le b_i \}$ (resp. $X \cap \{ f \le b_i^+ \}$) retracts by deformation to $X \cap \{f \le b_i \} \cap B_R^n(a)    $ (resp. $X \cap \{f \le b_i ^+\} \cap B_R^n(a)    $). Hence, we have:
$$\displaylines{
\quad \chi\big(\hbox{Lk}^\infty(X \cap \{ f \le b_i^+ \})\big) - \chi \big( \hbox{Lk}^\infty(X \cap \{ f \le b_i\}) \big)= \hfill \cr
\hfill \chi \big(X \cap \{ f \le b_i^+ \} \cap S_R^{n-1}    (a)\big) -\chi \big(X \cap \{ f \le b_i \} \cap S_R^{n-1}    (a)\big)= \quad
}$$
$$ \sum_{j :f(q_j^R) \in ]b_i,b_i^+[} \hbox{ind}(f,X \cap S_R^{n-1}    (a),q_j^R) =\sum_{j : \mu_j^R < 0 \atop f(q_j^R) \rightarrow b_i} \hbox{ind}(f,X \cap S_R^{n-1}    (a),q_j^R).$$
It remains to express $\sum_{j : \mu_j^R <0\atop f(q_j^R) \rightarrow -\infty}  \hbox{ind}(f,X \cap S_R^{n-1}    (a),q_j^R)$. Let us choose $R \gg 1$ and $b_0^+$ in $]-\infty,b_1[$ such that:
$$\left[  \cup_{j : f \rightarrow -\infty \atop along \ \mathcal{B}_j}  \mathcal{B}_j \right] \cap \left\{ \Vert x-a \Vert \ge R \right\} \subset f^{-1}(]-\infty,b_0^+[) \cap \left\{ \Vert x-a \Vert \ge R \right\} ,$$
$X \cap \{ f \le b_0^+ \}$ retracts by deformation to $X \cap \{f \le b_0^+ \} \cap B_R^n(a)    $ and:
$$\left[  \cup_{j : f \nrightarrow -\infty \atop along \ \mathcal{B}_j} \mathcal{B}_j \right] \cap \left\{ \Vert x-a \Vert \ge R \right\} \subset f^{-1}(]b_0^+,+\infty[) \cap \left\{ \Vert x-a \Vert \ge R \right\} .$$
By Theorem 3.1, we can write:
$$\displaylines{
\quad \chi \big(\hbox{Lk}^\infty(X \cap \{f \le b_0^+ \}) \big)= \chi \big(X \cap \{ f \le b_0^+ \} \cap S_R^{n-1}    (a)\} \big)= \hfill \cr
\hfil \sum_{j :f(q_j^R) \le b_0^+} \hbox{ind}(f,X \cap S_R^{n-1}    (a), q_j^R) =
\sum_{j : \mu_j^R <0\atop f(q_j^R) \rightarrow -\infty}\hbox{ind}(f,X \cap S_R^{n-1}    (a)    ,q_j^R). \quad \cr
}$$
To get the final result, we just remark that if $b \notin \Lambda_f^\le$ then:
$$\chi \big( \hbox{Lk}^\infty(X \cap \{f \le b^+ \})\big)-\chi \big( \hbox{Lk}^\infty(X \cap \{f \le b \}) \big)=0.$$ 
$\hfill \Box$

Similarly, $\Lambda_f^\ge =\{c_1,\ldots,c_s\}$ with $c_1<c_2< \cdots < c_s$ and:
$$\mathbb{R} \setminus \Lambda_f^\ge = ]-\infty,c_1 [ \  \cup \ ]c_1,c_2 [\  \cup \cdots \cup \  ]c_{s_1},c_s [ \ \cup \ ]c_s,+\infty [.$$
For each $i \in \{0,\ldots,s \}$, let $c_i^+$ be an element in $]c_i,c_{i+1}[$ with $c_0=-\infty$ and $c_{s+1}=+\infty$. 
\begin{theorem}
We have:
$$\displaylines{
\quad \chi(X)=\sum_{i=1}^l \hbox{\em ind}(-f,X,p_i) +\sum_{j=0}^s\chi \big(\hbox{\em Lk}^\infty(X \cap \{ f \ge c_j^+ \}) \big)  \hfill \cr
\hfill - \sum_{j=1}^s \chi \big(\hbox{\em Lk}^\infty(X \cap \{ f \ge c_j \}) \big). \quad \cr
}$$
\end{theorem} 
{\it Proof.} Same proof as Theorem 3.16. $\hfill \Box$

Let us write $\Lambda_f^= =\{d_1,\ldots,d_t \}$ with $d_1< d_2 < \ldots < d_t$ and:
$$\mathbb{R} \setminus \Lambda_f^= = ]-\infty,d_1 [ \ \cup \  ]d_1,d_2 [ \ \cup \cdots \cup \  ]d_{t_1},d_t [ \ \cup \ ]d_t,+\infty [.$$
For each $i \in \{0,\ldots,t \}$, let $d_i^+$ be an element in $]d_i,d_{i+1}[$.
\begin{corollary}
We have:
$$\displaylines{
\quad 2\chi(X) -\chi \big( \hbox{\em Lk}^\infty(X) \big)=\sum_{i=1}^l  \hbox{\em ind}(f,X,p_i)+\sum_{i=1}^l  \hbox{\em ind}(-f,X,p_i) + \hfill \cr
\hfill \sum_{j=0}^t\chi \big(\hbox{\em Lk}^\infty(X \cap \{ f = d_j^+ \}) \big) - \sum_{j=1}^t \chi \big( \hbox{\em Lk}^\infty(X \cap \{ f = d_j \}) \big). \quad \cr
}$$
\end{corollary}
{\it Proof.} Assume that  $\Lambda_f^\le \cup \Lambda_f^\ge = \Lambda_f^=$. We have:
$$\displaylines{
\quad \chi(X)=\sum_{i=1}^l \hbox{ind}(f,X,p_i) +\sum_{j=0}^t \chi \big( \hbox{Lk}^\infty(X \cap \{ f \le d_j^+ \}) \big) \hfill \cr
\hfill - \sum_{j=1}^t \chi \big(\hbox{Lk}^\infty(X \cap \{ f \le d_j \}) \big), \quad \cr
}$$
$$\displaylines{
\quad \chi(X)=\sum_{i=1}^l \hbox{ind}(-f,X,p_i) +\sum_{j=0}^t\chi \big( \hbox{Lk}^\infty(X \cap \{ f \ge d_j^+ \}) \big)  \hfill \cr
\hfill - \sum_{j=1}^t \chi \big( \hbox{Lk}^\infty(X \cap \{ f \ge d_j \}) \big). \quad \cr
}$$
Adding these two equalities and using the Mayer-Vietoris sequence, we obtain the result when $\Lambda_f^\le \cup \Lambda_f^\ge = \Lambda_f^=$ .
But if $d_j \notin \Lambda_f^=$ then $\chi \big(\hbox{Lk}^\infty(X \cap \{f= d_j^+ \}) \big)- \chi \big(\hbox{Lk}^\infty(X \cap \{f= d_j \}) \big)=0$. $\hfill \Box$

By Hardt's theorem, we know that there is a finite subset $B(f)$ of $\mathbb{R}$ such that $f_{\vert X \cap f^{-1}(B(f))}$ is a semi-algebraic locally trivial fibration. Hence outside $B(f)$, the function $\beta \mapsto \chi(X \cap \{f=\beta \})$ is locally constant. In the sequel, we will give formulas relating the topology of $X$ and the variations of topology in the fibres of $f$. Let us set $\tilde{B}(f) = f(\{p_1,\ldots,p_l\}) \cup \Lambda_f^\le \cup \Lambda_f^\ge$. This set is clearly finite. 
\begin{proposition}
If $\alpha \notin \tilde{B}(f)$ then the following functions:
$$\beta \mapsto \chi(X \cap \{f * \beta\}), \ * \in \{\le,=,\ge \},$$
are constant in a neighborhood of $\alpha$.
\end{proposition}
{\it Proof.}  We study the local behaviors of the numbers $\lambda_{f,\alpha}$ and $\mu_{f,\alpha}$, $\alpha \in \mathbb{R}$. We denote by $\alpha^+$ (resp. $\alpha^-$) an element of $]\alpha,+\infty[$ (resp. $]-\infty,\alpha[$) close to $\alpha$. If $R\gg 1$ is big enough and $\alpha^+$ is close enough to $\alpha$ then in $S_R^{n-1}    (a)     \cap \{\alpha \le f \le \alpha^+\}$, there is no points $q_j^R$ such that $\nabla (f_{\vert S})(q_j^R)=\mu_j^R \nabla ({\rho_a}_{\vert S} )(q_j^R)$ with $\mu_j^R>0$ because $f(q_j^R)$ decreases to $\alpha$ as $R$ tends to infinity. Hence if $\alpha^+$ is close enough to $\alpha$ then $\mu_{f,\alpha^+}=\mu_{f,\alpha}$.  In the same way, we can show that $\lambda_{f,\alpha^-}=\lambda_{f,\alpha}$. Applying this argument to $-f$ and $-\alpha$, we see that $\lambda_{-f,-\alpha^+}=\lambda_{-f,-\alpha}$ and $\mu_{-f,-\alpha^-}=\mu_{-f,-\alpha}$. If $\alpha \notin \tilde{B}(f)$ then, by Proposition 3.8, $\lambda_{f,\alpha^+}=\lambda_{f,\alpha}$, $\mu_{f,\alpha^-}=\mu_{f,\alpha}$, $\lambda_{-f,-\alpha^-}=\lambda_{f,-\alpha}$ and $\mu_{-f,-\alpha^+}=\mu_{-f,-\alpha}$. The formulas established in Proposition 3.6 and Corollary 3.7 enable us to conclude. $\hfill \Box$ 

Let us write $\tilde{B}(f)=\{\gamma_1,\ldots,\gamma_u\}$ and:
$$\mathbb{R} \setminus \tilde{B}(f)=]-\infty,\gamma_1[ \ \cup \ ]\gamma_1,\gamma_2[ \ \cup \cdots \cup \ ]\gamma_{u-1},\gamma_u[ \  \cup\ ]\gamma_u,+\infty[.$$
For $i \in \{0,\ldots,u\}$, let $\gamma_i^+$ be an element  of $]\gamma_i,\gamma_{i+1}[$ where $\gamma_0=-\infty$ and $\gamma_{u+1}=+\infty$. 
\begin{theorem}
We have:
$$\displaylines{
\quad \chi(X)=\sum_{i=1}^l  \hbox{\em ind}(f,X,p_i) +\sum_{i=1}^l \hbox{\em ind}(-f,X,p_i) +  \hfill \cr
\hfill \sum_{k=0}^u \chi \big(X \cap \{ f= \gamma_k^+\} \big) - \sum_{k=1}^u \chi \big( X \cap \{ f= \gamma_k\} \big) . \quad \cr
}$$
\end{theorem}
{\it Proof.} Let us assume first that $\tilde{B}(f)=f(\{p_1,\ldots,p_l\}) \cup \Lambda_f$, i.e that $\Lambda_f^\le \cup \Lambda_f^\ge =\Lambda_f$. 
For $k \in \{1,\ldots,u\}$, we have by Corollary 3.7:
$$\displaylines{
\quad \chi \big(X \cap \{f =\gamma_k \} \big)=\chi(X) - \sum_{i: f(p_i) > \gamma_k} \hbox{ind}(f,X,p_i) -  \sum_{i :f(p_i) < \gamma_k} \hbox{ind}(-f,X,p_i)  - \hfill \cr
 \hfill \lambda_{f,\gamma_k} -\lambda_{-f,-\gamma_k}, \quad \cr
 }$$
$$\displaylines{
\quad \chi \big(X \cap \{f =\gamma_k^+ \} \big)=\chi(X) -\sum_{i : f(p_i) > \gamma_k^+} \hbox{ind}(f,X,p_i) -  \sum_{i : f(p_i) < \gamma_k^+} \hbox{ind}(-f,X,p_i)   -\hfill \cr
\hfill \lambda_{f,\gamma_k^+} -\lambda_{-f,-\gamma_k^+}, \quad \cr
}$$
hence:
$$\displaylines{
\qquad \chi\big(X  \cap \{f =\gamma_k^+ \}\big)-\chi \big(X  \cap \{f =\gamma_k \} \big) = \hfill \cr
\hfill  - \sum_{i : f(p_i)=\gamma_k} \hbox{ind}(-f,X,p_i) -(\lambda_{f,\gamma_k^+} -\lambda_{f,\gamma_k}), \qquad \cr
}$$
because as already  noticed, $\lambda_{-f,-\gamma_k^+}=\lambda_{-f,-\gamma_k}$. 
If $\gamma_k$ does not belong to $\Lambda_f$ then $\lambda_{f,\gamma_k^+}=\lambda_{f,\gamma_k}$. If $\gamma_k$ belongs to $\Lambda_f$ then:
$$ \lambda_{f,\gamma_k} -\lambda_{f,\gamma_k^+}= \sum_{j :\mu_j^R <0 \atop f(q_j^R) \rightarrow \gamma_k}\hbox{ind}(f,X \cap S_R^{n-1}    (a) ,q_j^R).$$
Therefore,
$$\displaylines{
\quad \sum_{k=1}^u \chi \big(X \cap \{f=\gamma_k^+\} \big) -\chi \big(X \cap \{ f =\gamma_k \}\big)=-\sum_{i =1}^l \hbox{ind}(-f,X,p_i)+ \hfill \cr
\hfill \sum_{k:\gamma_k \in \Lambda_f} \sum_{j :\mu_j^R <0 \atop f(q_j^R) \rightarrow \gamma_k} \hbox{ind}(f,X \cap S_R^{n-1}    (a) ,q_j^R). \quad  \cr
}$$
By Corollary 3.7, we have:
$$\chi \big( X \cap \{f=\gamma_0^+\} \big)=\chi(X) -\sum_{i=1}^l \hbox{ind}(f,X,p_i) -\lambda_{f,\gamma_0^+} -\lambda_{-f,-\gamma_0^+}.$$
But we remark that:
$$\lambda_{f,\gamma_0^+} = \sum_{k: \gamma_k \in \Lambda_f} \sum_{j : \mu_j^R <0 \atop f(q_j^R) \rightarrow \gamma_k} \hbox{ind}(f,X \cap S_R^{n-1}    (a)    ,q_j^R),$$
because if $\nabla( f_{\vert S}) (q_j^R)= \mu_j^R \nabla ({\rho_a}_{\vert S}) (q_j^R) $ with $\mu_j^R <0$ then $f(q_j^R)> \gamma_0^+$ for $f(q_j^R)$ decreases to one of the $\gamma_i$'s. 
Similarly we see that $\lambda_{-f,-\gamma_0^+}=0$. Combining these equalities, we get the result when $\tilde{B}_f=\Lambda_f \cup f(\{p_1,\ldots,p_l\}) $.  But if $\gamma \notin \Lambda_f^\le \cup \Lambda_f^\ge \cup f(\Sigma_{f_{\vert X}})$ then $\chi \big(X \cap \{f =\gamma^+\} \big)-\chi \big(X \cap \{f=\gamma \} \big)=0$. 

$\hfill \Box$

\begin{theorem}
We have:
$$\chi(X)=\sum_{i=1}^l \hbox{\em ind}(f,X,p_i) + \sum_{k=0}^u \chi \big(X \cap \{ f \le \gamma_k^+ \} \big) -\sum_{k=1}^u \chi \big(X \cap \{ f \le \gamma_k \} \big),$$
$$\chi(X)=\sum_{i=1}^l  \hbox{\em ind}(-f,X,p_i) + \sum_{k=0}^u \chi \big(X \cap \{ f \ge \gamma_k^+ \} \big) -\sum_{k=1}^u \chi \big(X \cap \{ f \ge \gamma_k \} \big).$$
\end{theorem}
{\it Proof.} We prove the first equality in the case $\tilde{B}(f)=f(\{p_1,\ldots,p_l\})  \cup \Lambda_f$. For $k \in \{1,\ldots,u \}$, we have by Proposition 3.6:
$$\chi \big(X \cap \{ f \le \gamma_k \} \big) -\chi \big(X \cap \{ f = \gamma_k \} \big)= \sum_{i : f(p_i)< \gamma_k } \hbox{ind}(-f,X,p_i) + \lambda_{-f,-\gamma_k},$$
$$\chi \big(X \cap \{ f \le \gamma_k^+ \} \big) -\chi \big(X \cap \{ f = \gamma_k^+ \} \big)= \sum_{i : f(p_i)< \gamma_k^+ } \hbox{ind}(-f,X,p_i) + \lambda_{-f,-\gamma_k^+}.$$
Hence,
$$\displaylines{
\quad  \left[ \chi \big(X \cap \{ f \le \gamma_k^+ \} \big) -  \chi \big(X \cap \{ f \le \gamma_k \} \big) \right] - \hfill \cr
\hfill \left[ \chi \big(X \cap \{ f = \gamma_k^+ \} \big) -  \chi \big(X \cap \{ f = \gamma_k \} \big) \right] = \sum_{i:f(p_i)=\gamma_k} \hbox{ind}(-f,X,p_i), \quad \cr
}$$
so:
$$\displaylines{
\quad \sum_{k=1}^u \left[ \chi \big(X \cap \{ f \le \gamma_k^+ \} \big) -  \chi \big(X \cap \{ f \le \gamma_k \} \big) \right] - \hfill \cr
\hfill \sum_{k=1}^u  \left[ \chi \big(X \cap \{ f = \gamma_k^+ \} \big) -  \chi \big(X \cap \{ f = \gamma_k \} \big) \right] = \sum_{i=1}^l \hbox{ind}(-f,X,p_i). \quad \cr
}$$
But, we remark that $\chi \big(X \cap \{f \le \gamma_0^+ \} \big) =\chi \big(X \cap \{ f=\gamma_0^+\} \big)$ because $\lambda_{-f,-\gamma_0^+}=0$, which implies that:
$$\displaylines{
\quad \sum_{k=0}^u  \chi \big(X \cap \{ f \le \gamma_k^+ \} \big) -\sum_{k=1}^u  \chi \big(X \cap \{ f \le \gamma_k \} \big) = \hfill \cr
\quad \quad \sum_{k=0}^u  \chi \big(X \cap \{ f = \gamma_k^+ \} \big) -\sum_{k=1}^u  \chi \big(X \cap \{ f = \gamma_k \} \big)  + \sum_{i=1}^l \hbox{ind}(-f,X,p_i)= \hfill \cr
\hfill  \chi(X)-\sum_{i=1}^l \hbox{ind}(f,X,p_i). \quad \cr
}$$
$\hfill \Box$

\begin{corollary}
We have:
$$\displaylines{
\quad \sum_{k=0}^u  \chi \big(X \cap \{ f \ge \gamma_k^+ \} \big) - \chi \big(X \cap \{ f \le \gamma_k^+ \} \big) \hfill \cr
\quad \quad -\sum_{k=1}^u  \chi \big(X \cap \{ f \ge \gamma_k \} \big) - \chi \big(X \cap \{ f \le \gamma_k \} \big)=  \hfill \cr
\hfill \sum_{i=1}^l  \hbox{\em ind}(f,X,p_i)-\hbox{\em ind}(-f,X,p_i). \quad \cr
}$$
\end{corollary}
$\hfill \Box$ 

\section{Case $X=\mathbb{R}^n$}
In this section, we apply our previous results to the case $X=\mathbb{R}^n$. In this case $\hbox{ind}(f,X,p_i)=(-1)^n \hbox{ind}(-f,X,p_i)=\hbox{deg}_{p_i} \nabla f $, the local degree of $\nabla f$ at $p_i$ and $\sum_{i=1}^l \hbox{ind}(f,X,p_i)=\hbox{deg}_\infty \nabla f$, the degree of $\nabla f$ at infinity, i.e the topological degree of $\frac{\nabla f}{\vert \nabla f \vert} : S_R^{n-1}     \rightarrow S^{n-1}$ where $S_R^{n-1}    $ is a sphere of big radius $R$. Furthermore, $\mu_{f,\alpha}=(-1)^{n-1} \lambda_{-f,-\alpha}$ and $\mu_{-f,-\alpha}=(-1)^{n-1} \lambda_{f,\alpha}$. We can restate our result in this setting.
\begin{proposition}
For all $\alpha \in \mathbb{R}$, we have:
$$\chi\big(\{f \ge \alpha \} \big)-\chi \big(\{ f = \alpha \}\big)= \sum_{i :f(p_i)> \alpha } \hbox{\em deg}_{p_i} \nabla f + \lambda_{f,\alpha},$$
$$\chi \big(\{f \le \alpha \} \big)-\chi \big(\{ f = \alpha \} \big)= (-1)^n \sum_{i :f(p_i)< \alpha } \hbox{\em deg}_{p_i} \nabla f +(-1)^{n-1} \mu_{f,\alpha}.$$
\end{proposition}
$\hfill \Box$ 

\begin{corollary}
If $n$ is even then for all $\alpha \in \mathbb{R}$, we have:
$$\chi \big(\{f=\alpha \}\big)=1-\sum_{i : f(p_i) \not= \alpha } \hbox{\em deg}_{p_i} \nabla f -\lambda_{f,\alpha}+\mu_{f,\alpha},$$
$$\chi \big(\{f \ge \alpha \} \big)-\chi\big(\{ f \le \alpha \}\big)= \sum_{i :f(p_i)> \alpha } \hbox{\em deg}_{p_i} \nabla f -\sum_{i : f(p_i)<\alpha } \hbox{\em deg}_{p_i} \nabla f  +\lambda_{f,\alpha}+\mu_{f,\alpha}.$$
If $n$ is odd then for all $\alpha \in \mathbb{R}$, we have:
$$\chi \big(\{f=\alpha \}\big)=1- \sum_{i : f(p_i)> \alpha } \hbox{\em deg}_{p_i} \nabla f +\sum_{i : f(p_i)<\alpha } \hbox{\em deg}_{p_i} \nabla f -\lambda_{f,\alpha}+\mu_{f,\alpha},$$
$$\chi \big(\{f \ge \alpha \} \big)-\chi \big(\{ f \le \alpha \} \big)= \sum_{i : f(p_i) \not= \alpha } \hbox{\em deg}_{p_i} \nabla f + \lambda_{f,\alpha}+\mu_{f,\alpha}.$$
\end{corollary}
$\hfill \Box$ 

The above formulas can be viewed as real versions of results on the topology of the fibres of a complex polynomial (see for instance [Pa] or [ST]).
\begin{proposition}
If $n$ is even then, for all $\alpha \in \mathbb{R}$, we have:
$$\chi \big(\hbox{\em Lk}^\infty(\{f \le \alpha \}) \big)=\chi \big(\hbox{\em Lk}^\infty(\{f \ge \alpha \}) \big)=1-\hbox{\em deg}_\infty \nabla f -\lambda_{f,\alpha} + \mu_{f,\alpha},$$
$$\chi \big(\hbox{\em Lk}^\infty(\{f = \alpha \}) \big)=2-2\hbox{\em deg}_\infty \nabla f -2\lambda_{f,\alpha} +2 \mu_{f,\alpha}.$$
If $n$ is odd then, for all $\alpha \in \mathbb{R}$, we have:
$$\chi \big(\hbox{\em Lk}^\infty(\{f \le \alpha \})\big)=1-\hbox{\em deg}_\infty \nabla f -\lambda_{f,\alpha} + \mu_{f,\alpha},$$
$$\chi \big(\hbox{\em Lk}^\infty(\{f \ge \alpha \}) \big)=1+ \hbox{\em deg}_\infty \nabla f +\lambda_{f,\alpha} - \mu_{f,\alpha}.$$
\end{proposition}
$\hfill \Box$ 

We also obtain generalizations of Sekalski's formula [Se]. We keep the notations introduced in the general case.
\begin{theorem}
We have:
$$1= \hbox{\em deg}_\infty \nabla f + \sum_{i=0}^r \chi \big(\hbox{\em Lk}^\infty (\{ f \le b_i^+ \})\big) -\sum_{i=1}^r  \chi (\hbox{\em Lk}^\infty (\{ f \le b_i  \}) \big)=$$
$$(-1)^n \hbox{\em deg}_\infty \nabla f + \sum_{i=0}^s \chi \big( \hbox{\em Lk}^\infty (\{ f \ge c_i^+ \}) \big) -\sum_{i=1}^s \chi \big( \hbox{\em Lk}^\infty (\{ f \ge c_i  \}) \big).$$
If $n$ is even then we have:
$$2=2 \hbox{\em deg}_\infty \nabla f + \sum_{i=0}^t \chi \big(\hbox{\em Lk}^\infty (\{ f = d_i^+ \})\big) -\sum_{i=1}^t \chi \big(\hbox{\em Lk}^\infty (\{ f = d_i  \}) \big).$$
\end{theorem}
$\hfill \Box$ 

\begin{theorem}
If $n$ is even, we have:
$$1=2 \hbox{\em deg}_\infty \nabla f + \sum_{k=0}^u \chi \big(\{f = \gamma_k^+\} \big) -\sum_{k=1}^u \chi \big(\{f = \gamma_k\}\big),$$
$$1=\hbox{\em deg}_\infty \nabla f + \sum_{k=0}^u \chi \big(\{f \le \gamma_k^+\} \big) -\sum_{k=1}^u \chi \big(\{f \le \gamma_k\} \big),$$
$$1=\hbox{\em deg}_\infty \nabla f + \sum_{k=0}^u \chi \big(\{f \ge \gamma_k^+\} \big) -\sum_{k=1}^u \chi \big(\{f \ge \gamma_k\} \big),$$
$$  \sum_{k=0}^u \chi \big(\{f \ge \gamma_k^+\}\big)  -\chi \big(\{f \le \gamma_k^+\} \big)   =
\sum_{k=1}^u  \chi \big(\{f \ge \gamma_k \} \big) - \chi \big( \{f \le \gamma_k\} \big). $$
If $n$ is odd, we have:
$$1= \sum_{k=0}^u \chi \big(\{f = \gamma_k^+\} \big) -\sum_{k=1}^u \chi \big(\{f = \gamma_k\} \big),$$
$$1=\hbox{\em deg}_\infty \nabla f + \sum_{k=0}^u \chi \big(\{f \le \gamma_k^+\} \big) -\sum_{k=1}^u \chi \big(\{f \le \gamma_k\} \big),$$
$$1=-\hbox{\em deg}_\infty \nabla f + \sum_{k=0}^u \chi \big(\{f \ge \gamma_k^+\} \big) -\sum_{k=1}^u \chi \big(\{f \ge \gamma_k\} \big),$$
$$\displaylines{
\quad   \sum_{k=0}^u \chi \big(\{f \ge \gamma_k^+\} \big)  -\chi \big(\{f \le \gamma_k^+\} \big)   =
\sum_{k=1}^u  \chi \big(\{f \ge \gamma_k \} \big)-  \chi \big(\{f \le \gamma_k\} \big)  \hfill \cr
\hfill + 2 \hbox{\em deg}_\infty \nabla f. \quad \cr
}$$
\end{theorem}
$\hfill \Box$ 

\section{Application to generic linear functions}
We apply the results of Section 3 to the case of a generic linear function. 
Let $X \subset \mathbb{R}^n$ be a closed semi-algebraic set. For $v \in S^{n-1}$, let us denote by $v^*$ the function $v^*(x)= \langle v,x \rangle$. We are going to study the critical points of $v^*_{\vert X \cap S_R^{n-1}    (a)    }$ for $v$ generic and $R$ sufficiently big.

Let $\Gamma_1(X)$ be the subset of $S^{n-1}$ defined as follows: a  vector $v$ belongs to $\Gamma_1(X)$ if there exists a sequence $(x_k)_{k \in \mathbb{N}}$ such that $\Vert x_k \Vert \rightarrow +\infty$ and a sequence $(v_k)_{k \in \mathbb{N}}$ of vectors in $S^{n-1}$ such that $v_k \perp T_{x_k} S(x_k)$ and $v_k \rightarrow v$, where $S(x_k)$ is the stratum containing $x_k$. 

\begin{lemma}
The set $\Gamma_1(X)$ is a semi-algebraic set of $S^{n-1}$ of dimension strictly less than $n-1$.
\end{lemma}
{\it Proof.}  If we write 
$X =\sqcup_{\alpha \in A} S_\alpha$,
where $(S_\alpha)_{\alpha \in A}$ is a finite semi-algebraic Whitney stratification of $X$, then we see that $\Gamma_1(X)=
\sqcup_{\alpha \in A} \Gamma_1(S_\alpha)$. Hence it is enough to prove the lemma when $X$ is a smooth semi-algebraic manifold of dimension
$n-k$, $0< k< n$. 

Let us take $x=(x_1,\ldots,x_n)$ as a coordinate system for $\mathbb{R}^n$ and $(x_0,x)$ for
$\mathbb{R}^{n+1}$. Let $\varphi$ be the semi-algebraic diffeomorphism between $\mathbb{R}^n$ and $S^n \cap \{x_0 > 0\}$ given
by: 
$$\varphi (x) = \left( \frac{1}{\sqrt{1 + \Vert x \Vert^2}},\frac{x_1}{\sqrt{1 + \Vert x \Vert^2}}, \ldots,
\frac{x_n}{\sqrt{1 + \Vert x \Vert^2}} \right).$$
Observe that $(x_0,x)=\varphi(z)$ if and only if $z=\frac{x}{x_0}$. The set $\varphi(X)$ is a smooth
semi-algebraic set of dimension $n-k$. Let $M$ be the following semi-algebraic set:
$$M = \left\{ (x_0,x,y) \in \mathbb{R}^{n+1} \times \mathbb{R}^n \ \vert \ 
(x_0,x) \in \varphi(X) \hbox{ and } y \perp T_{\frac{x}{x_0}} X \right\}.$$
We will show that $M$ is a smooth manifold of dimension $n$. Let $p=(x_0,x,y)$ be a point in $M$ and let
$z=\varphi^{-1}(x_0,x)=\frac{x}{x_0}$. In a neighborhood of $z$, $X$ is defined by the vanishing of smooth functions
$g_1,\ldots,g_{k}$. For $i \in \{1,\ldots,k\}$, let $G_i$ be the smooth function defined by:
$$G_i(x_0,x)= g_i \left( \frac{x}{x_0} \right) = g_i (\varphi^{-1}(x_0,x)).$$
Then in a neighborhood of $(x_0,x)$, $\varphi(X)$ is defined by the vanishing of $G_1,\ldots,G_{k}$ and 
$x_0^2+x_1^2+\cdots+x_n^2-1 $. Note that for $i,k \in \{ 1,\ldots,n\}^2$, $\frac{\partial G_i}{\partial
x_k}(x_0,x)=\frac{1}{x_0} 
\frac{\partial g_i}{\partial x_k}(x).$
Hence in a neighborhood of $p$, $M$ is defined by the vanishing of $G_1,\ldots,G_{k}$,  
$x_0^2+x_1^2+\cdots+x_n^2-1 $ and the following minors $m_{i_1i_2\ldots i_{k+1}}$, $(i_1,\ldots,i_{k+1}) \in
\{1,\ldots,n\}^{k+1}$, given by:
$$m_{i_1i_2\ldots i_{k+1}}(x_0,x,y)= \left\vert \begin{array}{ccc}
\frac{\partial G_1}{\partial x_{i_1}}(x_0,x) & \cdots & \frac{\partial G_1}{\partial x_{i_{k+1}}}(x_0,x) \cr
\vdots & \ddots & \vdots \cr
\frac{\partial G_k}{\partial x_{i_1}}(x_0,x) & \cdots & \frac{\partial G_k}{\partial x_{i_{k+1}}}(x_0,x) \cr
y_{i_1} & \cdots & y_{i_{k+1}} \cr
\end{array} \right\vert.$$
Since rank$(\nabla
g_1,\ldots,\nabla g_k)=k$ at $z=\varphi^{-1}(x_0,x)$, one can assume that: 
$$\left\vert \begin{array}{ccc}
\frac{\partial G_1}{\partial x_{1}}(x_0,x) & \cdots & \frac{\partial G_1}{\partial x_{k}}(x_0,x) \cr
\vdots & \ddots & \vdots \cr
\frac{\partial G_k}{\partial x_{1}}(x_0,x) & \cdots & \frac{\partial G_k}{\partial x_{k }}(x_0,x) \cr
\end{array} \right\vert \not= 0 .$$
This implies that around $p$, $M$ is defined by the vanishing of $G_1,\ldots,G_k$, $m_{1\ldots k k+1},\ldots,m_{1\ldots k
n}$ and $x_0^2+x_1^2+\cdots+x_n^2-1$ (a similar argument is given and proved in [Dut1,\S 5]). It is straightforward to see that the gradient vectors of these functions are
linearly independent. Then $\bar{M}\setminus M$ is a semi-algebraic set of dimension less than $n$. If $\pi_y :
\mathbb{R}^{n+1} \times \mathbb{R}^n \rightarrow \mathbb{R}^n$ denotes the projection on the last $n$ coordinates, then we
have $\Gamma_1(X)= S^{n-1} \cap \pi_y (\bar{M}\setminus M)$. $\hfill \Box$ 

\begin{corollary}
Let $v$ be vector in $S^{n-1}$ and let $a \in \mathbb{R}^n$. If there exists a sequence $(x_k)_{k \in \mathbb{N}}$ of points in $X$ such that: 
\begin{itemize}
\item $\Vert x_k \Vert \rightarrow +\infty$,
\item $v \in N_{x_k} S(x_k) \oplus \mathbb{R}(x_k-a)$,
\item $\lim_{k \rightarrow + \infty} \vert v^*(x_k) \vert < + \infty$,
\end{itemize}
then $v$ belongs to $\Gamma_1(X)$ (Here $N_{x_k} S(x_k)$ is the normal space to the stratum $S(x_k)$).
\end{corollary}
{\it Proof.} We can assume that $v=e_1=(1,0,\ldots,0)$. In this case, $v^*=x_1$. Since the stratification is finite, we can 
assume that $(x_k)_{k \in \mathbb{N}}$ is a sequence of points lying in a stratum $S$. By the Curve
Selection Lemma  at infinity, there exists an analytic curve $p(t) : ]0,\varepsilon [ \rightarrow S$ such that $\lim_{t
\rightarrow 0} \Vert p(t) \Vert = + \infty$, $\lim_{t \rightarrow 0} p_1(t) < +\infty$ and for $t\in
]0,\varepsilon[$, $e_1$ belongs to the space $N_{p(t)} S \oplus \mathbb{R} ( p(t)-a)$. Let us consider the expansions as
Laurent series of the $p_i$'s:
$$p_i(t)=h_i t^{\alpha_i}+\cdots, \ i=1,\ldots,n.$$
Let $\alpha$ be the minimum of the $\alpha_i$'s. Necessarily, $\alpha <0$ and $\alpha_1 \ge 0$. It is straightforward to see
that $\Vert p(t) -a \Vert$ has an expansion of the form:
$$\Vert p(t) -a\Vert =b t^{\alpha}+\cdots, \  b>0.$$
Let us denote by $\pi_t$ the orthogonal projection onto $T_{p(t)} S$. For every $t \in ]0,\varepsilon [$, there exists a real
number $\lambda(t)$ such that:
$$\pi_t(e_1)= \lambda(t) \pi_t (p(t)-a )=\lambda(t) \Vert \pi_t (p(t)-a) \Vert \frac{\pi_t (p(t)-a)}{\Vert \pi_t (p(t)-a) \Vert}.$$
Observe that if $t$ is small enough, we can assume that $\pi_t(p(t)-a)$ does not vanish because $S_{\Vert p(t)-a \Vert}(a)$ intersects $S$ transversally. 
Using the fact that $p'(t)$ is tangent to $S$ at $p(t)$, we find that:
$$p_1'(t)=\langle p'(t),e_1 \rangle =\langle p'(t),\pi_t(e_1) \rangle= \lambda(t) \langle p'(t),p(t)-a\rangle.$$
This implies that ord$(\lambda )\ge \alpha_1 - 2\alpha$. Let $\beta$ be the order of $\Vert
\pi_t(p-a) \Vert$. Since $\Vert p(t) -a \Vert \ge \Vert \pi_t (p(t)-a )\Vert$, $\beta$
is greater or equal to $\alpha$. Finally we obtain that ord$(\lambda \Vert
\pi_t(p(t)-a) \Vert)$ is greater or equal to $\alpha_1 -2\alpha +\beta$, which is strictly positive. This proves the lemma. $\hfill \Box$ 

\begin{lemma}
There exists a semi-algebraic set $\Gamma_2(X) \subset S^{n-1}$ of dimension strictly less than $n-1$ such that if $v \notin \Gamma_2(X)$, then $v^*_{\vert X}$ has a finite number of critical points.
\end{lemma}
{\it Proof.} It is enough to prove the lemma for a semi-algebraic stratum $S$ of dimension $s<n$. Let $N_S$ be the following semi-algebraic set:
$$N_S = \{ (x,y) \in \mathbb{R}^n \times \mathbb{R}^n \ \vert \ x \in S \hbox{ and } y \perp T_x S \}.$$
Using the same kind of arguments as in Lemmas 2.2, 2.3, and 5.1, we see that $N_S$ is a smooth semi-algebraic manifold of dimension $n$. Let 
$$\begin{array}{ccccc}
\pi_y & : & N_S & \rightarrow  & \mathbb{R}^n \cr
   &  &  (x,y) & \mapsto & y  \cr
\end{array}$$
be the projection onto the last $n$ coordinates. The Bertini-Sard theorem implies that the set $D_{\pi_y}$ of critical values of $\pi_y$ is semi-algebraic of dimension strictly less than $n$. We take $\Gamma_2(X)= S^{n-1} \cap D_{\pi_y}$. $\hfill \Box$

Let us set $\Gamma(X)= \Gamma_1(X) \cup \Gamma_2(X)$, it is a semi-algebraic set of $S^{n-1}$ of dimension strictly less than $n-1$. If $v \notin \Gamma(X)$ then $v^*_{\vert X}$ admits a finite number of critical points $p_1,\ldots,p_l$. Moreover, if there is a family of points $q_j^R$ in $S \cap S_R^{n-1}    (a)    $ such that $\nabla v^*_{\vert S} (q_j^R)=\mu_j^R \nabla {\rho_a}_{\vert S}(q_j^R)$, where $S$ is  a stratum of $X$, and $\mu_j^R<0$ then $v^*(q_j^R) \rightarrow -\infty$ because $v \notin \Gamma_1$. Similarly if $\mu_j^R >0$ then $v^*(q_j^R) \rightarrow + \infty$. 
We conclude that the set $\Lambda_{v^*}$ is empty and that for all $\alpha \in \mathbb{R}$:
$$\lambda_{v^*,\alpha}= \mu_{v^*,\alpha}=\lambda_{-v^*,-\alpha}=\mu_{-v^*,-\alpha}=0.$$
Hence, we can restate the results of Section 3 in this setting and get relations between the topology of $X$ and the topology of generic hyperplane sections of $X$ (see [Ti2] for similar relations in the complex setting).
\begin{proposition}
If $v \notin \Gamma(X)$ then for all $\alpha \in \mathbb{R}$, we have:
$$\chi \big(X \cap \{v^* \ge \alpha\} \big) -\chi \big(X \cap \{v^* = \alpha\} \big)= \sum_{i : v^*(p_i)> \alpha} \hbox{\em ind}(v^*,X,p_i),$$
$$\chi \big(X \cap \{v^* \le \alpha\} \big) -\chi \big(X \cap \{v^* = \alpha\} \big)= \sum_{i : v^*(p_i)< \alpha} \hbox{\em ind}(-v^*,X,p_i),$$
$$\chi \big(X \cap \{v^* = \alpha\} \big)= \chi(X) - \sum_{i: v^*(p_i)> \alpha} \hbox{\em ind}(v^*,X,p_i) - \sum_{i : v^*(p_i)< \alpha}  \hbox{\em ind}(-v^*,X,p_i) ,$$
$$\displaylines{
\quad \chi \big(X \cap \{v^* \ge \alpha\} \big) -\chi \big(X \cap \{v^* \le \alpha\} \big)=  \hfill \cr
\hfill \sum_{i :v^*(p_i)> \alpha} \hbox{\em ind}(v^*,X,p_i)- \sum_{i : v^*(p_i)< \alpha}  \hbox{\em ind}(-v^*,X,p_i) . \hfill \cr
}$$
\end{proposition}
$\hfill \Box$ 

\begin{proposition}
If $v \notin \Gamma(X)$ then for all $\alpha \in \mathbb{R}$, we have:
$$\chi \big(\hbox{\em Lk}^\infty(X \cap \{v^* \le \alpha\}) \big) = \chi(X)-  \sum_{i =1}^l \hbox{\em ind}(v^*,X,p_i),$$
$$\chi \big(\hbox{\em Lk}^\infty(X \cap \{v^* \ge \alpha\}) \big) = \chi(X)-  \sum_{i =1}^l  \hbox{\em ind}(-v^*,X,p_i),$$
$$\displaylines{
\quad \chi \big( \hbox{\em Lk}^\infty(X \cap \{v^*=\alpha\}) \big) = 2\chi(X)-\chi \big(\hbox{\em Lk}^\infty(X) \big)\hfill \cr
\hfill -  \sum_{i =1}^l \hbox{\em ind}(v^*,X,p_i)-  \sum_{i =1}^l \hbox{\em ind}(-v^*,X,p_i). \quad \cr
}$$
\end{proposition}
$\hfill \Box$ 

Note that the functions $\beta \mapsto \chi(\hbox{Lk}^\infty(X \cap \{v^* ? \beta \}))$, $? \in \{\le,=,\ge\}$, are constant on $\mathbb{R}$.  Theorem 3.20, Theorem 3.21 and Corollary 3.22 are also valid in this context. They have the same formulation as in the general case with the difference that $\tilde{B}(v^*)=v^*(\{p_1,\ldots,p_l\}) $.

As an application, we will give a short proof of the Gauss-Bonnet formula for closed semi-algebraic sets. Let $\Lambda_0(X,-)$ be the Gauss-Bonnet measure on $X$ defined by:
$$\Lambda_0(X,U)=\frac{1}{\hbox{Vol}(S^{n-1})} \int_{S^{n-1}} \sum_{x \in U} \hbox{ind}(v^*,X,x) dv,$$
where $U$ is a Borel set of $X$. Note that if $x$ is not a critical point of $v^*_{\vert X}$ then ind$(v^*,X,x)=0$ and therefore that for $v \notin \Gamma(X)$, the sum 
$\sum_{x \in U} \hbox{ind}(v^*,X,x) $ is finite. The Gauss-Bonnet theorem for compact semi-algebraic sets is due to Fu [Fu] and Broecker and Kuppe [BK].
\begin{theorem}
If $X$ is a compact semi-algebraic set then:
$$\Lambda_0(X,X)=\chi(X).$$
\end{theorem}
Now assume that $X$ is just closed. Let $(K_R)_{R >0}$ be an exhaustive family of compact
Borel sets of $X$, that is a family $(K_R)_{R >0}$ of compact Borel sets of $X$ such that $\cup_{R >0} K_R
=X$ and $K_R \subseteqq K_{R'}$ if $R \le R'$. For every $R >0$, we have:
$$\Lambda_0 (X,X \cap K_R)= \frac{1}{\hbox{Vol} (S^{n-1})}\int_{S^{n-1}} \sum_{x \in X \cap K_R} \hbox{ind}(v^*,X,x) dv.$$
Moreover the following limit: $$\lim_{R \rightarrow +\infty} \sum_{x \in X \cap K_R}  \hbox{ind}(v^*,X,x),$$ is equal to $\sum_{x \in X}
\hbox{ind}(v^*,X,x)$ which is uniformly bounded by Hardt's theorem. Applying Lebesgue's theorem, we obtain:
$$\displaylines{
\quad \lim_{R \rightarrow + \infty} \Lambda_0 (X,X \cap K_R) = 
\frac{1}{\hbox{Vol} (S^{n-1})}\int_{S^{n-1}} \lim_{R \rightarrow + \infty} \sum_{x \in X \cap K_R} \hbox{ind}(v^*,X,x) dv=\hfill
\cr
\hfill \frac{1}{\hbox{Vol} (S^{n-1})}\int_{S^{n-1}} \sum_{x \in X} \hbox{ind}(v^*,X,x) dv. \quad \cr
}$$

\begin{definition}
We set: $$\Lambda_0 (X,X) = \lim_{R \rightarrow + \infty}  \Lambda_0(X,X \cap K_R),$$ where $(K_R)_{R >0}$ is an
exhaustive family of compact Borel sets of $X$.
\end{definition}

\begin{theorem}
If $X$ is a closed semi-algebraic set then:
$$\displaylines{
\qquad  \Lambda_0(X,X)= \chi(X)-\frac{1}{2} \chi\big(\hbox{\em Lk}^\infty(X)\big) \hfill \cr
\hfill -\frac{1}{2 \hbox{\em Vol}(S^{n-1})} 
\int_{S^{n-1}} \chi \big(\hbox{\em Lk}^\infty(X \cap \{v^*=0\}) \big)dv. \qquad \cr
}$$
\end{theorem}
{\it Proof.} We have:
$$\Lambda_0(X,X)=  \frac{1}{\hbox{Vol} (S^{n-1})}\int_{S^{n-1}} \sum_{x \in X} \hbox{ind}(v^*,X,x) dv =$$
$$ \frac{1}{\hbox{2Vol} (S^{n-1})}\int_{S^{n-1}} \sum_{x \in X} \hbox{ind}(v^*,X,x) + \hbox{ind}(-v^*,X,x) dv =$$
 $$\frac{1}{\hbox{2Vol} (S^{n-1})}\int_{S^{n-1}}  2\chi(X)-\chi\big(\hbox{Lk}^\infty(X)\big)- \chi \big(\hbox{Lk}^\infty(X \cap \{v^*=0 \})\big) dv ,$$
 by Proposition 5.5.  $\hfill \Box$

\end{document}